\documentclass[12pt]{article} 

\newcommand{\var}{\mathsf{var}}
\newcommand{\occ}{\mathsf{occ}}
\newcommand{\NAE}{\mathsf{NAE}}
\newcommand{\MU}{\mathsf{MU}}
\newcommand{\SMU}{\mathsf{SMU}}
\newcommand{\ROR}{\mathsf{ROR}}

\makeatletter
\def\mathcenterto#1#2{\mathclap{\phantom{#1}\mathclap{#2}}\phantom{#1}}
\let\old@widehat\widehat
\def\widehatto#1#2{\mathcenterto{#2}{\old@widehat{\mathcenterto{#1}{#2\,}}}}
\let\old@widehat\widehat
\def\widehatto#1#2{\mathcenterto{#2}{\old@widehat{\mathcenterto{#1}{#2\,}}}}
\makeatother

\makeatletter
\def\mathcenterto#1#2{\mathclap{\phantom{#1}\mathclap{#2}}\phantom{#1}}
\let\old@widetilde\widetilde
\def\widetildeto#1#2{\mathcenterto{#2}{\old@widetilde{\mathcenterto{#1}{#2\,}}}}
\let\old@widetilde\widetilde
\def\widetildeto#1#2{\mathcenterto{#2}{\old@widetilde{\mathcenterto{#1}{#2\,}}}}
\makeatother

\newcommand{\whbA}{\hspace{1.2pt}\widehatto{\bx}{\hspace{-1.2pt}\bA\hspace{1.2pt}}\hspace{-1.2pt}}
\newcommand{\wtbA}{\hspace{1.2pt}\widetildeto{\bx}{\hspace{-1.2pt}\bA\hspace{1.2pt}}\hspace{-1.2pt}}
\newcommand{\depth}{\mathsf{depth}}
\newcommand{\unsat}{\mathsf{unsat}}

\usepackage{amscd}
\usepackage{amsmath}
\usepackage{amssymb}
\usepackage{amstext}
\usepackage{amsthm}
\usepackage{bbold}
\usepackage{bm}
\usepackage{booktabs}
\usepackage{color}
\usepackage{colonequals}
\usepackage{comment}

\usepackage{easybmat}
\usepackage{etex}
\usepackage{enumitem}
\usepackage{framed}
\usepackage{authblk}
\usepackage{xcolor}
\usepackage[dvips,letterpaper,margin=1in]{geometry}
\usepackage{graphicx}
\usepackage{listings}
\usepackage{longtable}
\usepackage{mathtools}
\usepackage{multirow}
\usepackage{rotating}
\usepackage{setspace}
\usepackage{siunitx}
\usepackage{tabu}
\usepackage{verbatim}

\definecolor{cblack}{rgb}{0,0,0}
\definecolor{cblue}{rgb}{0.121569,0.466667,0.705882}    
\definecolor{corange}{rgb}{1.000000,0.498039,0.054902}  
\definecolor{cgreen}{rgb}{0.172549,0.627451,0.172549}   
\definecolor{cred}{rgb}{0.839216,0.152941,0.156863}     
\definecolor{cpurple}{rgb}{0.580392,0.403922,0.741176}  
\definecolor{cbrown}{rgb}{0.549020,0.337255,0.294118}   
\definecolor{cpink}{rgb}{0.890196,0.466667,0.760784}
\definecolor{cgray}{rgb}{0.498039,0.498039,0.498039}
\definecolor{cgreen2}{rgb}{0.7372549019607844, 0.7411764705882353, 0.13333333333333333}

\usepackage{hyperref}
\hypersetup{
  linkcolor  = cblue,
  citecolor  = cgreen,
  urlcolor   = corange,
  colorlinks = true,
}

\newtheorem{theorem}{Theorem}[section]

\newtheorem{lemma}[theorem]{Lemma}

\newtheorem{problem}[theorem]{Problem}
\newtheorem{definition}[theorem]{Definition}
\newtheorem{example}[theorem]{Example}
\newtheorem{proposition}[theorem]{Proposition}

\newtheorem{conjecture}[theorem]{Conjecture}
\theoremstyle{plain} 
\newcommand{\thistheoremname}{}
\newtheorem*{genericthm}{\thistheoremname}

\newcommand{\Komlos}{Koml\'{o}s}

\renewcommand{\emptyset}{\varnothing}

\makeatletter
\def\moverlay{\mathpalette\mov@rlay}
\def\mov@rlay#1#2{\leavevmode\vtop{%
   \baselineskip\z@skip \lineskiplimit-\maxdimen
   \ialign{\hfil$\m@th#1##$\hfil\cr#2\crcr}}}
\newcommand{\charfusion}[3][\mathord]{
    #1{\ifx#1\mathop\vphantom{#2}\fi
        \mathpalette\mov@rlay{#2\cr#3}
      }
    \ifx#1\mathop\expandafter\displaylimits\fi}
\makeatother

\newcommand{\RR}{\mathbb{R}}

\DeclareSymbolFont{bbold}{U}{bbold}{m}{n}
\DeclareSymbolFontAlphabet{\mathbbold}{bbold}

\newcommand{\ba}{\bm a}

\newcommand{\bx}{\bm x}
\newcommand{\by}{\bm y}

\newcommand{\bA}{\bm A}

\newcommand{\bM}{\bm M}

\DeclareSymbolFont{sfoperators}{OT1}{cmss}{m}{n}
\DeclareSymbolFontAlphabet{\mathsf}{sfoperators}
\makeatletter
\renewcommand{\operator@font}{\mathgroup\symsfoperators}
\makeatother

\DeclareMathOperator{\Tr}{Tr}

\DeclareMathOperator{\sgn}{sgn}

\newcommand{\Ex}{\mathop{\mathbb{E}}}  

\newcommand{\into}{\hookrightarrow}

\newif\ifnotes
\notestrue

\usepackage{relsize}

\title{The discrepancy of unsatisfiable matrices and a lower bound for the \Komlos\ conjecture constant}
\date{November 4, 2021}
\author{Dmitriy Kunisky\thanks{Email: \texttt{dmitriy.kunisky@yale.edu}. This work was supported by ONR Award N00014-20-1-2335 and a Simons Investigator Award to Daniel Spielman.}}
\affil{Department of Computer Science, Yale University}

\begin{document}

\maketitle

\begin{abstract}
    We construct simple, explicit matrices with columns having unit $\ell^2$ norm and discrepancy approaching $1 + \sqrt{2} \approx 2.414$.
    This number gives a lower bound, the strongest known as far as we are aware, on the constant appearing in the \Komlos\ conjecture.
    The \emph{unsatisfiable matrices} giving this bound are built by scaling the entries of clause-variable matrices of certain unsatisfiable Boolean formulas.
    We show that, for a given formula, such a scaling maximizing a lower bound on the discrepancy may be computed with a convex second-order cone program.
    Using a dual certificate for this program, we show that our lower bound is optimal among those using unsatisfiable matrices built from formulas admitting read-once resolution proofs of unsatisfiability.
    We also conjecture that a generalization of this certificate shows that our bound is optimal among all bounds using unsatisfiable matrices.
\end{abstract}

\thispagestyle{empty}

\clearpage

\tableofcontents

\thispagestyle{empty}

\clearpage

\setcounter{page}{1}

\section{Introduction}

The \emph{discrepancy} of a matrix $\bA \in \RR^{m \times n}$ is
\begin{equation}
    \Delta(\bA) \colonequals \min_{\bx \in \{ \pm 1\}^n} \|\bA\bx\|_{\infty}.
\end{equation}
The well-known \emph{\Komlos\ conjecture} of discrepancy theory concerns the quantities
\begin{align}
  K(n) &\colonequals \max\left\{\Delta(\bA): \bA \in \RR^{n \times n}, \text{ all columns of } \bA \text{ have } \ell^2 \text{ norm} \leq 1\right\}, \\
  K &\colonequals \sup_{n \geq 1} K(n).
\end{align}
\Komlos' conjecture is then that $K < \infty$.
Previous work on the conjecture has focused almost exclusively on upper bounds on $K(n)$, in particular elaborating on Banaszczyk's seminal result that $K(n) = O(\sqrt{\log n})$, which remains the best known upper bound \cite{Banaszczyk-1998-BalancingVectorsGaussianMeasures,BDGL-2018-GramSchmidtWalk,DGLN-2019-TowardsConstructiveBanaszczyk,BDG-2019-AlgorithmKomlosBanaszczyk,ALS-2021-DiscrepancySelfBalancingWalk}.

We are aware of very few works on lower bounds on $K(n)$ or attempts to disprove the \Komlos\ conjecture.
Hajela in \cite{Hajela-1988-KomlosConjectureFalseEvidence} showed that high-discrepancy instances may be found with respect to any subexponential number of sign vectors, which perhaps suggests that the conjecture might be false; \cite{CS-2021-NormsSignedSumsKomlos} recently improved slightly on this result.
Nikolov in~\cite{Nikolov-2013-KomlosVectorColorings} considered lower bounds arising from a natural semidefinite programming relaxation of $\Delta(\bA)$, but proved that these lower bounds never exceed~1.
In private communications, some researchers have suggested that the conjecture might hold with $K = 2$; however, Daniel Spielman recently communicated that an extensive numerical search yielded an (unpublished) example with $n = 15$ and $\Delta(\bA) \in [2.005, 2.006]$. To the best of our knowledge this is the state of the art at the time of writing.

Our main result is the following new lower bound on $K$.
\begin{theorem}
    \label{thm:K-lower}
    $K(2^k) \geq 2^{-k/2} + \sum_{a = 1}^k 2^{-a/2}$.
    Thus, $K \geq \sum_{a = 1}^{\infty} 2^{-a/2} = 1 + \sqrt{2} \approx 2.414$.
\end{theorem}
\noindent
Our proof of Theorem~\ref{thm:K-lower} will construct simple, explicit high-discrepancy matrices $\bA$.
These will belong to a more general class of constructions that we describe in Section~\ref{sec:unsat}, which convert unsatisfiable Boolean formulas into \emph{unsatisfiable matrices} whose discrepancy admits a convenient lower bound.
After proving the Theorem in Section~\ref{sec:proof}, we will also show in Section~\ref{sec:optimality} that, over bounds based on the special class of unsatisfiable matrices built from formulas admitting read-once resolution proofs of unsatisfiability, our bound is optimal.
We moreover conjecture that our construction is optimal over \emph{all} bounds based on unsatisfiable matrices, and show that this conjecture is implied by a concrete combinatorial claim about general resolution proofs.

We hope these results will stimulate future work investigating unsatisfiability and near-unsatisfiability as mechanisms through which collections of vectors can have high discrepancy.
It is perhaps surprising that this is a fruitful direction at all, since much recent work following Banaszczyk's, including that cited above, has viewed discrepancy as a geometric problem.
If the hard cases of discrepancy problems arise rather from logical obstacles as considered here, then it might be useful instead to ask whether discrepancy problems, even continuous ones such as \Komlos' conjecture, may in some general way be related back to logical or combinatorial problems like satisfiability.

\subsection{Notation}

For $x, y \in \RR$ we write $x \wedge y \colonequals \min\{x, y\}$.
We write $\sgn(x) \in \{-1, 0, +1\}$ for the sign of $x$, with $\sgn(0) = 0$.
For a matrix $\bA$, $\sgn(\bA)$ denotes the sign applied entrywise.
We emphasize that we will never consider the transpose of a matrix in this paper; the notation $\bA^T$ in Section~\ref{sec:unsat-from-trees} will instead denote a matrix associated to a tree $T$.
When we discuss rooted binary trees, we will always assume without further comment that they are \emph{full}, that is, that every internal vertex has exactly two children.
For $F$ a Boolean formula in conjunctive normal form, a \emph{variable} is one of the formal Boolean symbols occurring in $F$, while a \emph{literal} is $z$ or $\neg z$ for any variable $z$. We write $C \in F$ to indicate that $C$ is a clause of $F$, write $\var(F)$ for the set of variables in $F$, and write $\partial z$ with $z \in \var(F)$ for the set of $C \in F$ that contain $z$ or $\neg z$.

\section{Unsatisfiable formulas and matrices}
\label{sec:unsat}

We introduce the following definition of a type of matrix $\bA$ where the discrepancy admits a simple lower bound.

\begin{definition}
    \label{def:unsat}
    We call $\bA$ \emph{unsatisfiable} if, for all $\bx \in \{\pm 1\}^n$, there exists some $i \in [m]$ and $s \in \{ \pm 1\}$ such that, for all $j \in [n]$, either $A_{ij} = 0$ or $\sgn(x_j) = s \cdot \sgn(A_{ij})$.
\end{definition}
\noindent
That is, if we view the sign pattern of the rows of $\bA$ as describing ``templates'' for sign vectors---non-zero entries forcing a certain sign and zero entries acting as ``wildcards''---then every sign vector matches one of these templates, up to a global sign flip.
What is ``unsatisfiable'' is the condition that a vector avoids all of these templates at once.

\begin{example}
    \label{ex:unsat-matrices}
    Both of the following matrices are unsatisfiable, as are any matrices formed by multiplying their entries by non-negative scalars:
    \begin{equation}
        \left[\begin{array}{rrrr}
                1 & 1 & 1 & 0 \\
                1 & 1 & -1 & 0 \\
                1 & -1 & 0 & 1 \\
                1 & -1 & 0 & -1
              \end{array}\right]\,\, , \,\, \left[\begin{array}{rrrr}
                1 & 0 & 1 & 1 \\
                0 & 1 & 1 & -1 \\
                1 & 1 & -1 & 0 \\
                1 & -1 & 0 & -1
              \end{array}\right].
      \end{equation}
      The first is the kind of ``tree-like'' unsatisfiable matrix that we will analyze in detail in the sequel, while the second is outside this class; we discuss another remarkable property both share following our third proof of Theorem~\ref{thm:K-lower} in Section~\ref{sec:proof}.
\end{example}

The following will be our main tool.
Define
\begin{equation}
    \delta(\bA) \colonequals \min_{i = 1}^m \sum_{j = 1}^n |A_{ij}|.
\end{equation}
\begin{proposition}
    If $\bA$ is unsatisfiable, then $\Delta(\bA) \geq \delta(\bA)$.
\end{proposition}
\noindent
This is just because, by definition of $\bA$ being unsatisfiable, for each $\bx \in \{\pm 1\}^n$, $\bA\bx$ has an entry equal to plus or minus the $\ell^1$ norm of some row.
Because of this bound, it seems appealing to search among unsatisfiable matrices $\bA$ for ones with large $\delta(\bA)$, which will therefore also have large discrepancy.

\subsection{Connection with Boolean formulas}
\label{sec:nae}

Where can we find unsatisfiable matrices?
As our terminology suggests, it is helpful to reframe their definition in logical terms of Boolean formulas.
Let $z_1, \dots, z_n$ be formal Boolean variables.
For formulas $f_i$, define the ``not-all-equal'' or NAE clause $\NAE(f_{1}, \dots, f_{k})$ to be false if and only if its inputs are either all true or all false:
\begin{equation}
    \NAE(f_1, \dots, f_k) \colonequals (f_1 \vee \cdots \vee f_k) \wedge (\neg f_1 \vee \cdots \vee \neg f_k).
\end{equation}
Note that the NAE clause is invariant under negation, $\NAE(f_1, \dots, f_k) = \NAE(\neg f_1, \dots, \neg f_k)$.
We call a \emph{not-all-equal conjunctive normal form (NAE-CNF) formula} any conjunction of NAE clauses where all inputs are the $z_i$ or their negations, such as:
\begin{equation}
    F = \NAE(\neg z_1, z_2, \neg z_5) \wedge \NAE(z_2, \neg z_3, z_4, \neg z_6) \wedge \NAE(z_1, z_4, z_5).
\end{equation}
Any NAE-CNF formula $F$ may be encoded in its \emph{clause-variable matrix} $\bA^F$ with one row per NAE clause and one column per variable, such that $A^F_{ij}$ is 1 if $z_j$ occurs in clause $i$, is $-1$ if $\neg z_j$ occurs in clause $i$, and is zero otherwise.
For example, for the $F$ given above,
\begin{equation}
    \bA^F = \left[\begin{array}{rrrrrr} -1 & \hphantom{-}1 & 0 & 0 & -1 & 0 \\ 0 & 1 & -1 & \hphantom{-}1 & 0 & -1 \\ 1 & 0 & 0 & 1 & 1 & 0 \end{array}\right].
\end{equation}
We write $\bA^F$ for the matrix of formula $F$, and $F_{\bA}$ for the formula of matrix $\bA$.
The following is immediate from the definition of an unsatisfiable matrix.

\begin{proposition}
    $\bA$ is unsatisfiable (in the sense of Definition~\ref{def:unsat}) if and only if $F_{\sgn(\bA)}$ is unsatisfiable (in the usual sense of a Boolean formula).
\end{proposition}

\subsection{Unsatisfiable CNF to unsatisfiable NAE-CNF}
\label{sec:cnf-to-nae}

We may also build unsatisfiable matrices almost as directly from ordinary unsatisfiable conjunctive normal form (CNF) formulas, rather than ones involving NAE clauses.
To see this, note that the CNF formula
\begin{equation}
    (f_{i_{1,1}} \vee \cdots \vee f_{i_{1, c_1}}) \wedge \cdots \wedge (f_{i_{m,1}} \vee \cdots \vee f_{i_{m, c_m}})
\end{equation}
is satisfiable if and only if the following NAE-CNF formula, equipped with an extra ``dummy'' variable $z^{\prime}$, is satisfiable:
\begin{equation}
    \NAE(z^{\prime}, f_{i_{1,1}}, \dots, f_{i_{1, c_1}}) \wedge \cdots \wedge \NAE(z^{\prime}, f_{i_{m,1}}, \dots, f_{i_{m, c_m}}).
\end{equation}
Thus, for an unsatisfiable CNF formula $F$, we write $\bA^F$ for the matrix formed by first converting it to the equivalent unsatisfiable NAE-CNF formula, and then building the associated matrix.
In fact, if we associate to $F$ itself a clause-variable matrix $\wtbA^F$ as we did to an NAE-CNF formula, by recording which variables occur positively or negatively in which clauses, then $\bA^F$ is formed simply by appending to $\wtbA^F$ a column of 1's (which corresponds to the variable $z^{\prime}$).
We adopt the convention that this concatenation happens on the left:
\begin{equation}
    \bA^F = \left[\begin{array}{c|c} 1\,\,\, & \\ 1\,\,\, & \\ \vdots\,\,\, & \hspace{1cm} \wtbA^F  \\ 1\,\,\, & \\ 1\,\,\, \end{array} \hspace{0.8cm}\right].
\end{equation}
We will mostly focus in the remainder of the paper on unsatisfiable matrices built from CNF formulas via this transformation; in Section~\ref{sec:normalization} we will briefly discuss the question of whether working with NAE-CNF formulas not of this form can improve on our results.

\subsection{Unsatisfiable formulas from trees}
\label{sec:unsat-from-trees}

We have reduced the task of building unsatisfiable matrices to that of producing unsatisfiable CNF formulas.
We introduce two classes of such formulas whose structure is especially well-understood.
The following definition is not crucial to our reasoning, but we present it to explain the names of these classes of formulas before giving the reformulations that will be more useful.
The reader may see, e.g., \cite{BHM-2009-HandbookSatisfiability} for a detailed treatment of these topics.
\begin{definition}
    A formula $F$ is \emph{minimally unsatisfiable (MU)} if $F$ is unsatisfiable, but removing any clause from $F$ yields a satisfiable formula.
    $F$ is furthermore \emph{strongly minimally unsatisfiable (SMU)} if adding any literal to any clause (not already containing that literal) yields a satisfiable formula.
\end{definition}
\noindent
The following quantity turns out to give a natural grading of unsatisfiable formulas.
\begin{definition}
  The \emph{deficiency} of an unsatisfiable formula $F$ is the difference between the number of clauses and the number of variables.
  The set of MU formulas $F$ with deficiency $\ell$ is denoted $\MU(\ell)$ and the set of SMU formulas formulas with deficiency $\ell$ is denoted $\SMU(\ell)$.
\end{definition}
\noindent
A folklore result, sometimes attributed to Tarsi, says that all unsatisfiable formulas have deficiency at least one.

\begin{figure}
    \begin{center}
    \begin{tabular}{m{7cm}m{1cm}m{4cm}}
      \includegraphics[scale=0.7]{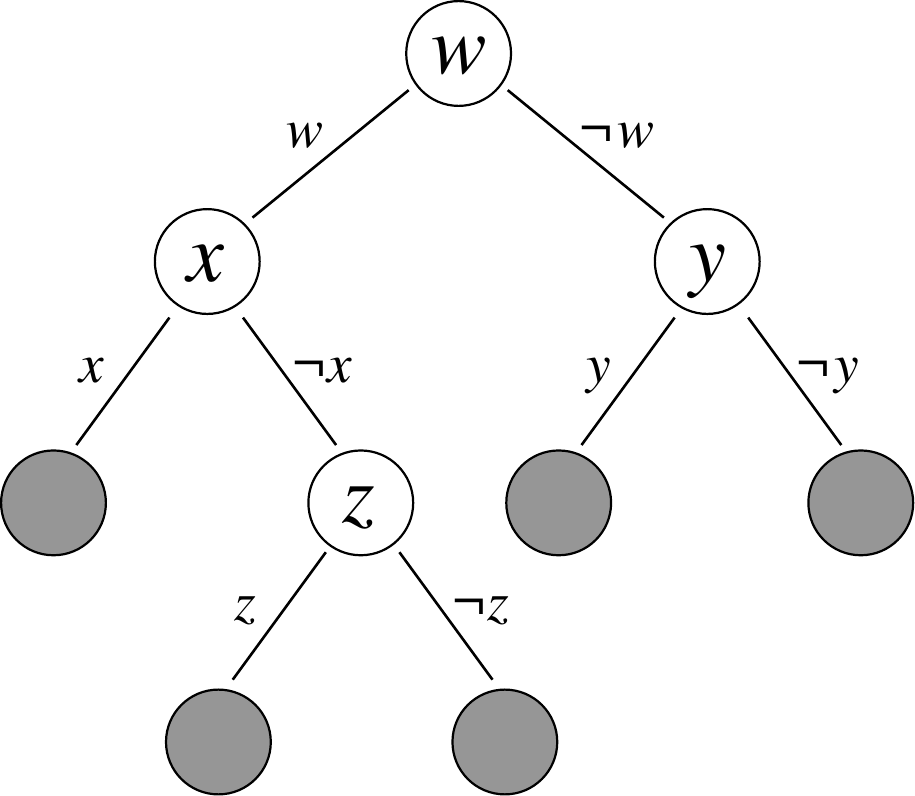} &  $\mathlarger{\mathlarger{\mathlarger{\mathlarger{\rightarrow}}}}$ & \vspace{-3.6em} \begin{tabular}{l}
                                                                          $(w \vee x) \, \, \wedge$ \\
                                                                          $(w \vee \neg x \vee z) \, \, \wedge$ \\
                                                                          $(w \vee \neg x \vee \neg z) \, \, \wedge$ \\
                                                                          $(\neg w \vee y) \, \, \wedge$ \\
                                                                          $(\neg w \vee \neg y)$
                                                                          \end{tabular}
    \end{tabular}
\end{center}
    \caption{We illustrate the tree construction of an $\SMU(1)$ formula. The clauses of the formula on the right, read from top to bottom, correspond to the leaves of the tree on the left, read from left to right.}
    \label{fig:tree-smu}
\end{figure}

The following is a complete and concrete description of the classes $\MU(1)$ and $\SMU(1)$.
While there are some structural results about the classes $\MU(\ell)$ and $\SMU(\ell)$ for $\ell \geq 2$, these results quickly become far weaker for growing $\ell$ (e.g., beyond $\ell \geq 5$; see \cite{BK-2009-MinimalUnsatisfiabilityAutarkies} for a survey).
\begin{definition}
  \label{def:tree-formula}
  Let $T$ be a rooted binary tree, with each internal vertex $v$ labelled by a distinct Boolean variable $z_v$, one leafward edge from $v$ labelled by $z_v$, and the other labelled by $\neg z_v$.
  Label each leaf of $T$ by a clause formed by the disjunction of the literals appearing as edge labels on the path joining that leaf to the root.
  We write $F_T$ for the CNF formula given by the conjunction of the clauses appearing on the leaves of $T$.
\end{definition}
\noindent
We illustrate this construction in Figure~\ref{fig:tree-smu}.
Perhaps surprisingly, this class of constructions exactly describes $\SMU(1)$, and a minor variation describes $\MU(1)$.

\begin{proposition}[Lemma C.5 of \cite{Kullmann-2000-MatroidTheorySAT}]
    \label{prop:mu-smu}
    For every $F \in \SMU(1)$, there exists a rooted binary tree $T$ as in Definition~\ref{def:tree-formula} such that $F = F_T$.
    For every $F \in \MU(1)$, there exists $F^{\prime} \in \SMU(1)$ such that $F$ is formed by removing occurrences of literals from clauses of $F^{\prime}$ such that every literal still occurs in some clause.
\end{proposition}
\noindent
We note that, by the requirement on the occurrence of literals in $F$, it is impossible to ``collapse'' two clauses and make them identical when we start with $F^{\prime} \in \SMU(1)$ and follow the construction above.
While \cite{Kullmann-2000-MatroidTheorySAT} gives a convenient summary, the result for $\SMU(1)$ originates with \cite{AL-1986-MinimalUnsatisfiableFormulas}, and essentially the same statement for $\MU(1)$, albeit with a more complicated proof and not making explicit reference to trees, appeared in \cite{DDB-1998-MinimalUnsatisfiabilityCNF}.

Let us briefly outline the simple and elegant proof that each $F_T$ is unsatisfiable: given an assignment of the variables, follow the path from the root given by, at each step, moving to the vertex that is labelled with a false value under this assignment.
Then, every literal in the disjunctive clause of the leaf thus reached is false, so the assignment is not satisfying.
This holds for all assignments, so $F_T$ is unsatisfiable.
The same argument also applies to any $\MU(1)$ formula formed by removing occurrences of literals from $F_T$.

This tree construction of unsatisfiable formulas has proved useful in other searches for extremal unsatisfiable CNF formulas.
Perhaps most notably, the line of work \cite{Gebauer-2012-NeighborhoodConjectureSAT,GST-2016-LocalLemmaTightSAT} designed trees $T$ such that associated $\MU(1)$ formulas make the greatest ``degree'' of any variable (the number of clauses containing it) as small as is possible, to leading order, over all unsatisfiable $k$-CNF formulas.
In this sense, the $F_T$ and their truncations appear to include the ``most regular'' possible unsatisfiable formulas.
This will be relevant for our purposes, since we too will want to find unsatisfiable formulas that are as regular as possible, albeit in a ``smoother'' sense we discuss in Section~\ref{sec:normalization}.

For the sake of brevity, we write $\bA^T \colonequals \bA^{F_T}$ for the unsatisfiable matrix associated to the formula $F_T \in \SMU(1)$.
Later we will consider optimizing this construction over multiplications of the entries of $\bA^T$ by arbitrary non-negative scalars, which implicitly allows the ``literal elimination'' forming $\MU(1)$ formulas from $\SMU(1)$ ones.

We note also that, since in any binary tree the number of internal vertices is one smaller than the number of leaves, the matrix $\bA^T$ constructed per Sections~\ref{sec:nae} and \ref{sec:cnf-to-nae} from $F_T$ will be square.
This is not important (a rectangular matrix could be padded with zeros to be square without changing its discrepancy, for instance), but coincides nicely with the traditional statement of the \Komlos\ conjecture over square matrices.

\section{Application to the \Komlos\ conjecture}
\label{sec:proof}

We will now prove Theorem~\ref{thm:K-lower}, by producing unsatisfiable $\whbA \in \RR^{n \times n}$ with $\delta(\whbA)$ large and columns of unit norm.
These will have $\sgn(\whbA) = \bA^T$ for complete binary trees $T$, and will have their entry magnitudes adjusted to satisfy the norm constraint.
Before proceeding, let us introduce one last piece of terminology: we call a matrix \emph{normalized} if all of its columns have $\ell^2$ norm at most 1 (recalling that this is the condition of the \Komlos\ conjecture).
As a notational convention, we write $\whbA$ for matrices guaranteed to be normalized (possibly with various superscripts) and $\bA$ for other matrices.

\subsection{Proof of Theorem~\ref{thm:K-lower}}

We give three proofs.
All three work with the following sequence of matrices.
Write $T_k$ for the complete binary tree of depth $k$.
The first few matrices $\bA^{T_k}$ produced by the procedure of Section~\ref{sec:unsat-from-trees} from these trees are
\begin{equation}
     \arraycolsep=4.0pt\def\arraystretch{1.0}
     \underbrace{\left[\begin{array}{c} 1 \end{array}\right]}_{\bA^{T_0}} \, \, , \, \,
     \underbrace{\left[\begin{array}{rr} 1 & 1 \\ 1 & -1 \end{array}\right]}_{\bA^{T_1}} \, \, , \, \,
     \underbrace{\left[\begin{array}{rrrr}
                        1 & 1 & 1 & 0 \\
                        1 & 1 & -1 & 0 \\
                        1 & -1 & 0 & 1 \\
                        1 & -1 & 0 & -1
                      \end{array}\right]}_{\bA^{T_2}} \, \, , \, \,
     \underbrace{\left[\begin{array}{rrrrrrrr}
                        1 & 1 & 1  & 0  & 1  & 0  & 0 & 0 \\
                        1 & 1 & 1  & 0  & -1 & 0  & 0 & 0 \\
                        1 & 1 & -1 & 0  & 0  & 1  & 0 & 0 \\
                        1 & 1 & -1 & 0  & 0  & -1 & 0 & 0 \\
                        1 & -1 & 0 & 1  & 0  & 0  & 1 & 0\\
                        1 & -1 & 0 & 1  & 0  & 0  & -1 & 0 \\
                        1 & -1 & 0 & -1 & 0  & 0  & 0 & 1\\
                        1 & -1 & 0 & -1 & 0  & 0  & 0 & -1
                      \end{array}\right]}_{\bA^{T_3}}.
            \end{equation}
(Here we view $F_{T_0}$ as the empty CNF formula.)
We recall that the leading column all of whose entries equal 1 corresponds to the dummy variable added to a CNF formula to get an equivalent NAE-CNF formula.

The pattern is likely intuitively clear, but let us give a more precise description.
The columns of $\bA^{T_k}$, besides the first column of the dummy variable, correspond to internal vertices of $T_k$, which we label with binary strings, starting with the empty string $\emptyset$ at the root and adding a 0 to get the string of the left child and a 1 to get that of the right child for each internal vertex.
The rows of $\bA^{T_k}$ correspond to the leaves of $T_k$, which we label by continuing the same procedure.
For a binary string $s$, let us write $|s|$ for the length, $s(i)$ for the $i$th character, and $s < t$ if $s$ is a strict prefix of the string $t$.
Then, it is simple to verify that the procedure of Section~\ref{sec:unsat-from-trees} gives
\begin{equation}
    A^{T_k}_{st} = \left\{\begin{array}{rl} 1 & \text{if } t < s \text{ and } s(|t| + 1) = 0, \\ -1 & \text{if }t < s \text{ and } s(|t| + 1) = 1, \\
                              0 & \text{otherwise}. \end{array}\right.
\end{equation}

Let $\whbA^{T_k}$ be formed by dividing every column of $\bA^{T_k}$ by its $\ell^2$ norm.
Counting the vertices at various depths in $T_k$ and the numbers of leaves in their subtrees gives that, for $k \geq 1$, there are two columns in $\bA^{T_k}$ whose norm is $2^{k/2}$ (one is the column of the dummy variable), and $2^a$ columns whose norm is $2^{(k - a)/2}$ for each $1 \leq a \leq k - 1$.
Also, every row has a non-zero entry in both columns of the first kind, and a non-zero entry in one column of each subsequent kind.
Thus the $\ell^1$ norm of each row of $\whbA^{T_k}$ is the same, and equals
\begin{equation}
    \delta(\whbA^{T_k}) = 2 \cdot 2^{-k/2} + \sum_{a = 1}^{k - 1} 2^{-(k - a)/2} = 2^{-k/2} + \sum_{a = 1}^k 2^{-a/2}.
\end{equation}
(For $k = 0$, the only column is that of the dummy variable, and the adjusted calculation gives the correct result $\delta(\whbA^{T_0}) = 1$ as well.)
Therefore, Theorem~\ref{thm:K-lower} is proved provided we can show that the $\bA^{T_k}$ are unsatisfiable.
It is this claim that we give several proofs of.

The first proof is immediate based on the background we have established.

\begin{proof}[Proof 1]
    The $\bA^{T_k}$ are the unsatisfiable matrices associated to the unsatisfiable formulas $F_{T_k}$ built from the trees $T_k$ per Section~\ref{sec:unsat-from-trees}.
\end{proof}

The matrices $\bA^{T_k}$ are also well-known in the signal processing literature as the transposes of those representing the discrete wavelet transform of the Haar wavelet (the normalized $\whbA^{T_k}$ give the orthogonal version of this transform).
Thus these matrices operate as a ``synthesis'' operator combining wavelet coefficients into a signal.
Our claim may then be seen as saying that perfectly ``flat'' collections of wavelet coefficients must give rise to signals that are at least somewhat ``spiky.''
The usual recursive treatment of these matrices in signal processing also affords a different proof of Theorem~\ref{thm:K-lower}.

\begin{proof}[Proof 2]
    Up to reordering the rows and columns, the normalized matrices satisfy the recursion
    \begin{equation}
        \whbA^{T_{k + 1}} = \frac{1}{\sqrt{2}}\left[\begin{array}{rr} \whbA^{T_k} & \bm I_{2^{k}} \\ \whbA^{T_k} & -\bm I_{2^k} \end{array}\right]
    \end{equation}
    for all $k \geq 1$.
    It is simple to show by induction that this coincides with the construction given earlier.
    To prove $\whbA^{T_k}$ is unsatisfiable, we prove by induction that $\whbA^{T_k}\bx$ has some entry equal exactly to $\delta(\whbA^{T_k})$ for all $\bx \in \{ \pm 1\}^{2^k}$.
    The statement clearly holds for $k = 0$.
    For $\bx_1, \bx_2 \in \{\pm 1\}^{2^{k}}$ we have
    \begin{equation}
        \whbA^{T_{k + 1}} \left[\begin{array}{c} \bx_1 \\ \bx_2 \end{array}\right] = \frac{1}{\sqrt{2}}\left[\begin{array}{c} \whbA^{T_k}\bx_1 + \bx_2 \\ \whbA^{T_k}\bx_1 - \bx_2 \end{array}\right].
    \end{equation}
    By the inductive hypothesis, some entry of $\whbA^{T_k}\bx_1$ equals $\pm \delta(\whbA^{T_k})$, say in coordinate $i$.
    Then one of coordinate $i$ or coordinate $2^{k} + i$ above equals $\pm (1 + \delta(\whbA^{T_k})) / \sqrt{2} = \pm \delta(\whbA^{T_{k + 1}})$, completing the induction.
\end{proof}
\noindent
The same argument more generally implies that the map
\begin{equation}
    \phi(\bA) \colonequals \frac{1}{\sqrt{2}}\left[\begin{array}{rr} \bA & \bm I \\ \bA & -\bm I \end{array}\right]
\end{equation}
has the property of preserving both normalization and unsatisfiability of matrices.
We have
\begin{equation}
    \delta(\phi(\bA)) = \frac{1 + \delta(\bA)}{\sqrt{2}},
\end{equation}
so in fact a recursive sequence $\bA^{(k)} = \phi(\bA^{(k - 1)})$ initialized from \emph{any} normalized unsatisfiable $\bA^{(0)}$ (even one with $\delta(\bA^{(0)}) > 1 + \sqrt{2}$) will have $\delta(\bA^{(k)}) \to 1 + \sqrt{2}$.

Our final proof of unsatisfiability takes advantage of a further special property of the $\bA^{T_k}$ beyond unsatisfiability: the sets of sign vectors that align with each of their rows (up to a sign flip) are \emph{disjoint}.

\begin{proof}[Proof 3]
    Let us say that $\bx \in \{\pm 1\}^{2^k}$ \emph{matches} row $i$ of $\bA^{T_k}$ if there exists $s \in \{\pm 1\}$ such that, for all $j$ with $A^{T_k}_{ij} \neq 0$, we have $x_j = s \cdot A^{T_k}_{ij}$.
    $\bA^{T_k}$ is unsatisfiable if and only if each $\bx$ matches some row of $\bA^{T_k}$.
    No $\bx$ can match two different rows of $\bA^{T_k}$, since any pair of rows have a pair of equal non-zero entries in the first coordinate and a pair of opposite non-zero entries in some coordinate.
    Also, since each row has $k + 1$ nonzero entries, each row matches $2^{2^k - k - 1} \cdot 2 = 2^{2^k} / 2^k$ different vectors of $\{\pm 1\}^{2^k}$.
    Since there are $2^k$ rows, every $\bx \in \{\pm 1\}^{2^k}$ matches exactly one row of $\bA^{T_k}$.
\end{proof}
\noindent
This proof points to a more geometric interpretation of unsatisfiability.
The negation of an unsatisfiable CNF formula, or CNF contradiction, is a \emph{disjunctive normal form (DNF)} formula that is true for every variable assignment, thus called a DNF tautology.
A DNF tautology corresponds to a covering of the hypercube $\{\pm 1\}^n$ by facets (the same as our description of Definition~\ref{def:unsat} in terms of ``templates'' for sign vectors).
Such coverings have been studied occasionally both from the geometric and the logical perspectives \cite{AL-1986-MinimalUnsatisfiableFormulas,FHK-1993-PackingLinesHypercube,SS-2000-DNFTautologiesVariableOccurrences,Szorenyi-2008-DNFTautologiesConflictBound}.

Our proof shows that the rows of $\bA^{T_k}$ describe parallel pairs of facets of the hypercube $\{ \pm 1\}^n$ that \emph{partition} the hypercube vertices rather than merely covering them.
Such DNF tautologies are called \emph{disjoint} in \cite{Szorenyi-2008-DNFTautologiesConflictBound}, and \cite{AL-1986-MinimalUnsatisfiableFormulas} observed that $\SMU(1)$ formulas always correspond to disjoint covers.
This property makes it easier to check unsatisfiability, requiring only that we verify disjointness (the rows' agreeing in at least one non-zero entry and disagreeing in at least one non-zero entry) and carrying out a counting argument.

We mention in passing a tempting attack on the \Komlos\ conjecture along these lines: the second $4 \times 4$ matrix in Example~\ref{ex:unsat-matrices} is one example of an unsatisfiable matrix with this partitioning property that does not arise from $\SMU(1)$, or for that matter from any CNF formula, and that has the further special property of having an equal number of non-zero entries in each column.
If such $2^k \times 2^k$ unsatisfiable matrices with $k + 1$ non-zero entries in each row \emph{and} each column existed for arbitrarily large $k$ then they would have discrepancy at least $\sqrt{k + 1} \gtrsim \sqrt{\log 2^k}$, disproving the \Komlos\ conjecture and showing that Banaszczyk's discrepancy bound is tight.
However, counting the requisite agreeing and disagreeing pairs of coordinates shows that such regular matrices cannot exist once $k > 2$.
It may seem appealing to look for nearly-regular matrices instead, but, based on Conjecture~\ref{conj:dual-cert} below, we believe that this approach will not disprove the \Komlos\ conjecture.

\section{Optimality over read-once resolution formulas}
\label{sec:optimality}

It is natural to ask whether our lower bound is optimal among those based on tree constructions of unsatisfiable formulas.
In Theorem~\ref{thm:complete-is-optimal} below we will show that this is indeed the case: any matrix $\whbA$ constructed by normalizing the entries of $\bA^F$ with $F \in \MU(1)$ to have columns of norm at most 1 with arbitrary entry magnitudes has $\delta(\whbA) < 1 + \sqrt{2}$.
In fact, more is true: the same holds for the larger class of all $F$ admitting \emph{read-once resolution} proofs of unsatisfiability.

Towards stating our result, we review some basic definitions concerning resolution.
Let $F$ be an unsatisfiable CNF formula.
It is well-known that $F$ admits a \emph{tree resolution} proof of unsatisfiability.
This is a rooted binary tree $T$ whose vertices are labelled by disjunctive clauses.
Each leaf is labelled with some clause of $F$ (with repetitions allowed) and the root is labelled with the empty clause.
Whenever the two children of an internal vertex are labelled with clauses $C$ and $D$, there is some variable $z$ such that it is possible to write $C = z \vee G$ and $D = \neg z \vee H$.
In this case, the internal vertex is labelled with $G \vee H$ and we say that the variable $z$ is \emph{resolved} at that vertex.
Under this rule, if an assignment makes the clauses labelling the two children of an internal vertex true, then it also makes the clause labelling the internal vertex true.
Since the root is labelled with the empty clause, which is never true, the existence of such a tree certifies the unsatisfiability of $F$.

If there exists such a resolution proof where each clause occurs as the label of at most one leaf of $T$, then we say that $F$ has a \emph{read-once resolution (ROR)} proof of unsatisfiability, and we write $F \in \ROR$.
The tree representation of an $\MU(1)$ formula from Definition~\ref{def:tree-formula} and Proposition~\ref{prop:mu-smu} in fact also gives a tree resolution proof of unsatisfiability of the same formula, where the label of each internal vertex is the variable that is resolved at that vertex.
Moreover, these are read-once proofs, so we have $\SMU(1) \subset \MU(1) \subset \ROR$, and each inclusion is strict; the former is easy to check, while a small example for the latter is given in \cite{BK-2009-MinimalUnsatisfiabilityAutarkies}.
In particular, as that example shows, it can be that a formula is not MU but is ROR and ceases to be ROR if any clauses are removed; that is, sometimes one may make an unsatisfiable formula ROR by adding some superfluous clauses and violating minimal unsatisfiability.

Our second main result is then as follows.
\begin{theorem}
    \label{thm:complete-is-optimal}
    Let $F \in \ROR$, and let $\whbA$ be any matrix formed by multiplying the entries of $\bA^{F}$ by arbitrary non-negative scalars so that the columns have $\ell^2$ norm at most~1.
    Then, $\delta(\whbA) < 1 + \sqrt{2}$.
\end{theorem}

\subsection{Optimal normalization of an unsatisfiable matrix}
\label{sec:normalization}

We first show that, for any unsatisfiable matrix $\bA$, there is an efficient convex optimization algorithm for maximizing $\delta(\whbA)$ over those matrices $\whbA$ that are normalized by multiplying the entries of $\bA$ by non-negative scalars.
That is, tuning the magnitudes of the entries $\bA$ is a fairly routine matter once their signs are fixed.\footnote{The study of those properties of a matrix that depend only on its sign pattern is sometimes called \emph{qualitative matrix analysis (QMA)}, which has played a role in other results in satisfiability; see \cite{BK-2009-MinimalUnsatisfiabilityAutarkies} for some discussion.}
This fact both will play a role in our proof of Theorem~\ref{thm:complete-is-optimal}, as we will use the convex duality of this optimization problem, and is a useful computational device (we have used it for searching for high-discrepancy unsatisfiable matrices and verifying Theorem~\ref{thm:complete-is-optimal} numerically), so we present it in generality before moving on to our argument.

\begin{lemma}
    \label{lem:opt-norm}
    Let $F$ be any unsatisfiable NAE-CNF formula on $n$ variables.
    Let $C_1, \dots, C_m \subseteq [n]$ contain the indices of variables appearing in each clause, and $V_1, \dots, V_n \subseteq [m]$ contain the indices of clauses that each variable appears in.
    Over all matrices formed by multiplying the entries of $\bA^F$ by non-negative scalars such that the columns have $\ell^2$ norm at most 1, let $\whbA^F$ be the matrix that maximizes $\delta(\whbA^F)$.
    Then, this value is given by either of the following dual optimization problems, which are convex second-order cone programs (SOCPs):
    \begin{align}
      \delta(\whbA^F)
        &= \left\{ \begin{array}{ll} \text{maximize} & \min_{i = 1}^m \sum_{j \in C_i} a_{ij} \\ \text{subject to} & \sqrt{\sum_{i \in V_j} a_{ij}^2} \leq 1 \text{ for all } j \in [n] \end{array}\right\} \\
        &= \left\{\begin{array}{ll} \text{minimize} & \sum_{j = 1}^n \sqrt{\sum_{i \in V_j} \mu_{i}^2} \\ \text{subject to} & \sum_{i = 1}^m \mu_{i} = 1 \end{array} \right\}.
    \end{align}
    The variables in the first (``primal'') program are $a_{ij} \in \RR$ indexed by $i \in [m]$ and $j \in [n]$ such that $j \in C_i$, while the variables in the second (``dual'') program are $\mu_i$ indexed only by $i \in [m]$.
    The constraints $a_{ij} \geq 0$ and $\mu_i \geq 0$ may be added to the primal and dual respectively without changing the value.
\end{lemma}
\noindent
The primal program is just a direct formulation of optimizing over possible normalizations of $\bA^F$.
Obtaining the dual is a routine convex duality computation; we give the details in Appendix~\ref{app:socp} by writing the above as an equivalent semidefinite program, computing its dual, and then converting that dual back to an SOCP.

Another way to see that the primal is a convex problem without working with second-order constraints is to set $b_{ij} \colonequals \sqrt{a_{ij}}$ and square the constraint, making this a maximization of the minimum of sums of square roots under linear constraints.
Since the square root function is concave, the minimum of sums of square root functions (of linear inputs) is concave also.
Thus the primal amounts to maximizing a concave function under linear constraints.

In light of the dual program in Lemma~\ref{lem:opt-norm}, the \Komlos\ conjecture for the smaller quantity $\delta(\cdot)$ rather than $\Delta(\cdot)$ and for the special case of unsatisfiable matrices reads as follows.
\begin{problem}
    \label{prob:unsat-opt}
    Does there exist $K_{\unsat} > 0$ such that, for any unsatisfiable NAE-CNF formula $F$, there exists a distribution of ``masses'' $\mu_C \geq 0$ over clauses with $\sum_{C \in F} \mu_C = 1$ and
    \begin{equation}
        \sum_{z \in \var(F)} \sqrt{\sum_{C \in \partial z} \mu_C^2} \leq K_{\unsat} \, ?
    \end{equation}
\end{problem}
\noindent
The optimal $K_{\unsat}$ is a lower bound on $K$ from \Komlos' conjecture, and the cases of CNF and NAE-CNF formulas are equivalent up to an adjustment of $K_{\unsat}$.
Namely, if we expand each NAE clause $C$ of an NAE-CNF formula to two disjunctive clauses $C_1$ and $C_2$, obtaining an unsatisfiable CNF formula, then for any clause weights $\bm \nu$ for the CNF formula we may set $\mu_C \colonequals \nu_{C_1} + \nu_{C_2}$, and since $\mu_C^2 \leq 2(\nu_{C_1}^2 + \nu_{C_2}^2)$ the objective value of clause weights $\bm \mu$ will be at most $\sqrt{2}$ times that of $\bm \nu$.
If Conjecture~\ref{conj:dual-cert} below holds, then Problem~\ref{prob:unsat-opt} admits a positive solution with $K_{\unsat} = \sqrt{2}(1 + \sqrt{2}) = 2 + \sqrt{2} \approx 3.414$; thus, in this case, unsatisfiable NAE-CNF formulas could at best prove the lower bound $K \geq 2 + \sqrt{2}$.
This would be an improvement on our lower bound $K \geq 1 + \sqrt{2}$, but numerical searches have shown no indication that such NAE-CNF formulas exist.

We note that the work of \cite{Gebauer-2012-NeighborhoodConjectureSAT,GST-2016-LocalLemmaTightSAT} mentioned earlier addressed a similar question: it was conjectured that, in any unsatisfiable $k$-CNF formula, there is some variable that occurs in at least $\Omega(2^k / k)$ clauses, i.e., some $z \in \var(F)$ with $|\partial z| = \Omega(2^k / k)$.
The results of \cite{Gebauer-2012-NeighborhoodConjectureSAT} constructed trees achieving this, and \cite{GST-2016-LocalLemmaTightSAT} later even achieved the optimal constant (see also \cite{KST-1993-kSATVariableOccurrence, SS-2000-DNFTautologiesVariableOccurrences, HS-2005-UnsatisfiableSATFewOccurrences, HS-2006-UnsatisfiableNCFFewOccurrences} for earlier progress and \cite{HKSS-2014-PositionalGames} for another presentation).
Problem~\ref{prob:unsat-opt} also asks how ``regular'' unsatisfiable CNF formulas can be, albeit in a ``smoother'' way.
To see that this is the general thrust of the question, note for example that if there are $m$ variables, $n$ clauses, and $\mu_C = 1 / n$ for all $C \in F$, then the left-hand side is the average of $\sqrt{|\partial z|}$ over variables $z$, divided by the clause density $n / m$.

\subsection{Proof of Theorem~\ref{thm:complete-is-optimal}}
\label{sec:complete-is-optimal}

We give two proofs, which share the following initial observations.
Without loss of generality, we may suppose $\whbA = \whbA^F$ is the optimal normalization of $\bA^F$ maximizing $\delta(\whbA^F)$.
The first step is to reduce to the case $F \in \MU(1)$, and then to the case $F \in \SMU(1)$.

Let $T$ be a read-once resolution tree for $F$, and view each internal vertex as labelled by the variable that is resolved at that vertex.
Consider repeatedly performing the following operation: if any internal vertex label is repeated, find an occurrence of a repeated label whose subtree contains no further occurrences of that label, and change that variable and all its occurrences in that subtree to a new variable added to $F$.
The resulting $T$ remains a resolution tree of the new formula.
On $\bA^F$, this operation has the effect of splitting one column into two, with each non-zero entry occurring in exactly one of the two columns in the same coordinate.
This only increases $\delta(\whbA^F)$, since any normalization of the old matrix remains a valid normalization of the new matrix (with the entries of the split column retaining their normalized magnitudes) with the same value of $\delta(\cdot)$.
Thus we may repeat this operation until no internal vertex labels in $T$ are repeated, only increasing $\delta(\whbA^F)$; at the end, all internal vertex labels will be distinct, and the resulting formula will be in $\MU(1)$.
Thus without loss of generality we may suppose $F \in \MU(1)$.

Next, note that, by Proposition~\ref{prop:mu-smu}, for any $F \in \MU(1)$ there exists $F^{\prime} \in \SMU(1)$ such that $\bA^F$ is produced by zeroing some entries of $\bA^{F^{\prime}}$.
Thus we may further without loss of generality suppose $F \in \SMU(1)$, i.e., $F = F_T$ and thus $\bA^F = \bA^T$ for some binary tree $T$.

We will prove Theorem~\ref{thm:complete-is-optimal} by producing a dual certificate $\bm \mu$ for the convex program in Lemma~\ref{lem:opt-norm}.
To motivate this dual certificate, it is helpful to identify a good dual certificate for the case of the complete binary tree $T = T_k$.
For $v$ a vertex of a given tree, let us write $\depth(v)$ for its depth, and for $T$ a tree let us write $L(T)$ for the set of leaves and $I(T)$ for the set of internal vertices.
We write $\ell \leq i$ if $\ell$ lies in the subtree rooted at $i$ (or ``below'' $i$).
In $T_k$, the optimal dual bound is achieved by the dual certificate $\bm\nu \in \RR^{2^k}$ where we set all $\nu_{\ell} \colonequals 2^{-k}$.
Indeed, with these choices the objective value is
\begin{align}
  \sqrt{\sum_{\ell \in L(T_k)} \nu_{\ell}^2} + \sum_{i \in I(T_k)} \sqrt{ \sum_{\ell \leq i} \nu_{\ell}^2 }
  &= 2^{-k/2} + \sum_{i \in I(T_k)}\sqrt{ 2^{k - \depth(i)} \cdot 2^{-2k}} \nonumber \\
  &= 2^{-k/2} + \sum_{i \in I(T_k)} 2^{-(k + \depth(i)) / 2} \label{eq:Tk-internal-ell2}\\
  &= 2^{-k/2} + \sum_{a = 0}^{k - 1} 2^a \cdot 2^{-(k + a) / 2} \nonumber \\
  &= 2^{-k/2} + \sum_{a = 1}^{k}2^{-a / 2}.
\end{align}

If $T$ has depth at most $k$, we view $T$ as a subtree of $T_k$, so that some of the leaves of $T$ correspond to internal vertices of $T_k$ where the larger tree has been ``trimmed.''
We choose $\mu_{\ell}$ for $\ell \in T$ so that the squared mass of $\ell$ is the proportional to the sum of the squared masses of $\ell^{\prime} \in T_k$ that lies under $\ell$:
\begin{align}
  \mu_{\ell}^{(0)} &\colonequals \sqrt{\sum_{\ell^{\prime} \in L(T_k): \ell^{\prime} \leq \ell} \nu_{\ell^{\prime}}^2} = 2^{-(k + \depth(\ell)) / 2}, \\
  \mu_{\ell} &\colonequals \frac{\mu_{\ell}^{(0)}}{\sum_{\ell \in L(T)} \mu_{\ell}^{(0)}},
\end{align}
where we borrow the first computation from \eqref{eq:Tk-internal-ell2}.

We will then have for any $i \in I(T)$ that
\begin{equation}
    \sqrt{\sum_{\ell \in L(T): \ell \leq i} \mu_{\ell}^2} = \frac{1}{\sum_{\ell \in L(T)} \mu_{\ell}^{(0)}} \sqrt{\sum_{\ell \in L(T): \ell \leq i} \sum_{\ell^{\prime} \in L(T_k): \ell^{\prime} \leq \ell} \nu_{\ell^{\prime}}^2} = \frac{2^{-(k + \depth(i)) / 2}}{\sum_{\ell \in L(T)} 2^{-(k + \depth(\ell)) / 2}},
\end{equation}
borrowing the same computation again.
The dual bound achieved by $\bm \mu$ is the sum of these terms, with that of the root vertex repeated:
\begin{equation}
    \label{eq:dual-bd}
    \delta(\whbA^T) \leq \frac{2^{-k/2} + \sum_{i \in I(T)} 2^{-(k + \depth(i)) / 2}}{\sum_{\ell \in L(T)} 2^{-(k + \depth(\ell)) / 2}} = \frac{1 + \sum_{i \in I(T)} 2^{-\depth(i) / 2}}{\sum_{\ell \in L(T)} 2^{-\depth(\ell) / 2}}.
\end{equation}
It thus suffices to show that the quantity on the right-hand side, now just a combinatorial quantity associated to a binary tree $T$, is always smaller than $1 + \sqrt{2}$.
Having used $T_k$ to motivate our construction, we now work only with $T$ and write $I = I(T)$ and $L = L(T)$ below.

\begin{proof}[Proof 1]
    Let us write $I_a$ and $L_a$ for the sets of internal vertices and leaves respectively at depth $a$ in $T$.
    We have $|I_0| = 1$, $|L_0| = 0$, and the $|L_a|$ may be computed from the $|I_a|$ as
    \begin{equation}
        |L_a| = 2|I_{a - 1}| - |I_a|
    \end{equation}
    for $a \geq 1$.
    Using this, we may rewrite the denominator of \eqref{eq:dual-bd} in terms of the numerator as
    \begin{align}
      \sum_{\ell \in L} 2^{-\depth(\ell) / 2}
      &= \sum_{a = 1}^{\infty} |L_a| 2^{-a/2} \nonumber \\
      &= \sqrt{2}\sum_{a = 0}^{\infty} |I_a| 2^{-a/2} - \sum_{a = 1}^{\infty} |I_a|2^{-a/2} \nonumber \\
      &= 1 + (\sqrt{2} - 1)\sum_{a = 0}^{\infty} |I_a| 2^{-a/2} \nonumber \\
      &= 1 + (\sqrt{2} - 1)\sum_{i \in I}2^{-\depth(i) / 2}. \label{eq:L-I-relation}
    \end{align}
    Thus the bound \eqref{eq:dual-bd} is, setting $X \colonequals \sum_{i \in I}2^{-\depth(i) / 2} > 0$,
    \begin{equation}
        \delta(\whbA^T) \leq \frac{1 + X}{1 + (\sqrt{2} - 1)X},
    \end{equation}
    which is a strictly increasing function of $X > 0$ whose limit as $X \to \infty$ is $1 + \sqrt{2}$.
\end{proof}

This proof is simple and direct, but suffers from the drawback that it seems quite difficult to extend to the generalized construction we propose later in Section~\ref{sec:conj-general}.
We therefore also offer the following, seemingly more elaborate proof, which however offers a technique that might be more readily generalized.
It also gives natural probabilistic interpretations of the unusual-looking quantities appearing in Theorem~\ref{thm:K-lower} as
\begin{align}
  2^{-k/2} + \sum_{a = 1}^k 2^{-a/2} &= \Ex_{X \sim \mathsf{Geom}(1/2)} \sqrt{2}^{X \wedge k}, \\
  1 + \sqrt{2} &= \Ex_{X \sim \mathsf{Geom}(1/2)} \sqrt{2}^{X},
\end{align}
where $\mathsf{Geom}(1/2)$ is the geometric distribution of the number of fair coin tosses made before a toss lands heads.

\begin{proof}[Proof 2]
    Let $f: I \into L$ be an injective map such that $f(i) \leq i$ for all $i \in I$.
    That is, $f$ maps each internal vertex to a distinct leaf lying below it.
    Write $p_f(i) = p_f(f(i))$ for the length of the path between $i$ and $f(i)$, with $p_f(\ell) = -\infty$ if $\ell$ is not in the image of $f$.
    For any such $f$,
    \begin{equation}
        \sum_{i \in I} 2^{-\depth(i) / 2} = \sum_{\ell \in L} 2^{-\depth(\ell) / 2} \cdot 2^{p_f(\ell) / 2}.
    \end{equation}
    Suppose further that $\rho$ is a probability distribution over such $f$.
    We then have
    \begin{equation}
        \sum_{i \in I} 2^{-\depth(i) / 2} = \sum_{\ell \in L} 2^{-\depth(\ell) / 2} \Ex_{f \sim \rho}\left[2^{p_f(\ell) / 2}\right].
    \end{equation}

    We choose the random $f \sim \rho$ sampled as follows.
    We first construct a partition of $T$ into edge-disjoint paths from internal vertices to leaves lying beneath those vertices.
    We start two paths from the root.
    At each non-root internal vertex, with probability $1/2$ each we continue the path arriving at that vertex from above to the left or to the right, and start a new path in the other direction.
    Next, we delete one of the two paths touching the root, uniformly at random.
    In the remaining set of paths, there is exactly one path starting at each internal vertex, and each path ends at a distinct leaf.
    Thus, for each $i \in I$, we let $f(i)$ be the leaf endpoint of the unique path whose internal vertex endpoint is $i$.
    We illustrate this construction in Figure~\ref{fig:path-partition}.

    \begin{figure}
        \begin{center}
            \includegraphics[scale=0.5]{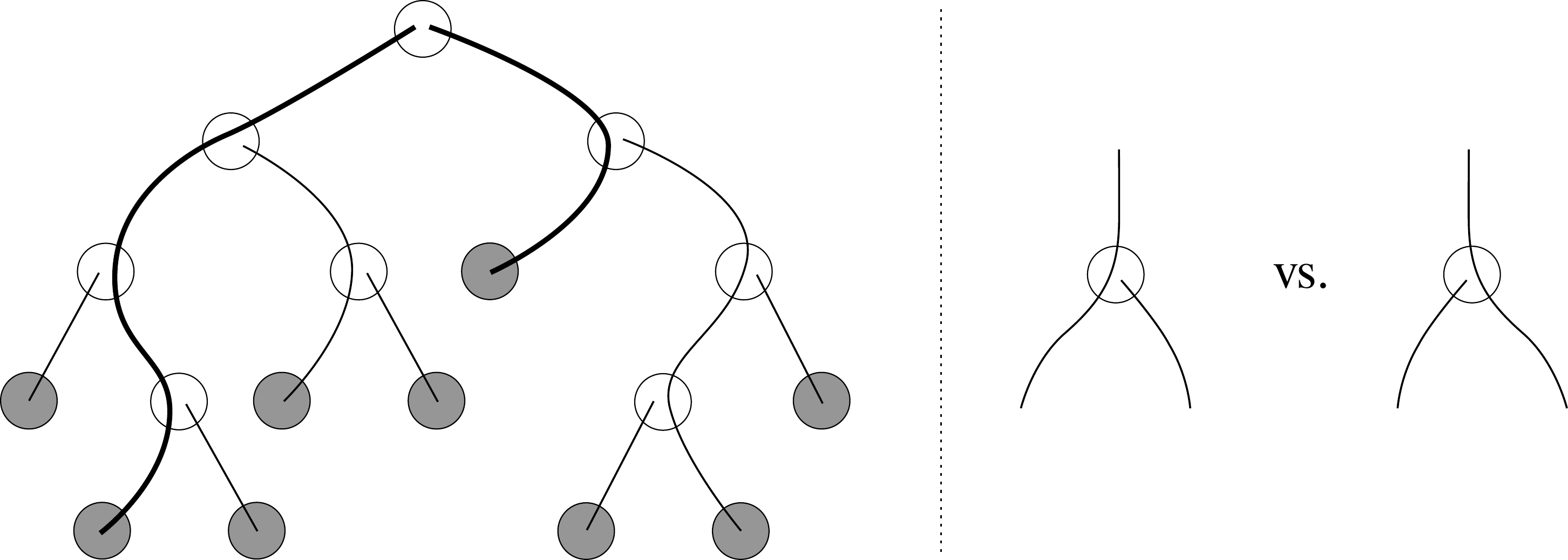}
        \end{center}
        
        \caption{In the left panel, we show a partition of a binary tree into edge-disjoint paths used in the second proof of Theorem~\ref{thm:complete-is-optimal}. The two bold paths are those attached to the root; to form the map $f$, one of the two will be discarded from the partition. In the right panel, we show the local choice made at each non-root internal vertex in constructing a random such partition.} 
        \label{fig:path-partition}
    \end{figure}
    
    For any $\ell$, the law of the random variable $p_f(\ell)$ is that of $(X \wedge \depth(\ell)) \cdot Y$, where $X \sim \mathsf{Geom}(1/2)$, $Y = 1$ whenever $X < \depth(\ell)$, and $Y$ is chosen uniformly at random from $\{-\infty, 1\}$ if $X = \depth(\ell)$.
    The case $Y = -\infty$ here corresponds to $\ell$ being the endpoint of one of the paths touching the root in the initial partition, but that path then being deleted when the mapping $f$ is determined, so that $\ell$ is not in the image of $f$.
    In particular, $p_f(\ell) = j$ with probability $2^{-j}$ for $j = 1, \dots, \depth(\ell)$, and $p_f(\ell) = -\infty$ with probability $2^{-\depth(\ell)}$.
    Thus,
    \begin{equation}
        \Ex_{f \sim \rho}\left[2^{p_f(\ell) / 2}\right] = \sum_{j = 1}^{\depth(\ell)} 2^{-j/2} = \frac{2^{-1/2} - 2^{-(\depth(\ell) + 1) / 2}}{1 - 2^{-1/2}} = (1 - 2^{-\depth(\ell) / 2})(1 + \sqrt{2}).
    \end{equation}
    
    Substituting and using that $\sum_{\ell \in L} 2^{-\depth(\ell)} = 1$, we find
    \begin{equation}
        \sum_{i \in I} 2^{-\depth(i) / 2} = (1 + \sqrt{2})\sum_{\ell \in L} (2^{-\depth(\ell) / 2} - 2^{-\depth(\ell)}) = (1 + \sqrt{2})\left(\sum_{\ell \in L} 2^{-\depth(\ell) / 2} - 1\right),
    \end{equation}
    which is just a rewritten version of the main relation \eqref{eq:L-I-relation} from the first proof, and the rest of the proof is completed as before.
\end{proof}

\subsection{The stick game}

Though we do not use it in our arguments, we mention a pleasant reinterpretation of searching for high-discrepancy matrices built by normalizing $\bA^T$ for binary trees $T$.
The following describes the simultaneous construction both of the tree and the normalization.

Consider a game of woodland solitaire, where the player is equipped with an endless supply of sticks of unit length.
The game is played on several piles of sticks, starting with one empty pile.
At every step, the player can choose one pile, add a new stick of unit length to it, break every stick in that pile into two sticks of arbitrary length, and place one from each pair into a new pile.
For instance, representing piles by tuples, the game may proceed like so:
\begin{equation}
    \underline{\bigg(\,\bigg)} \to \underline{\bigg(\frac{1}{3}\bigg)} \bigg(\frac{2}{3}\bigg) \to \bigg(\frac{1}{12}, \frac{1}{2}\bigg) \underline{\bigg(\frac{1}{4}, \frac{1}{2}\bigg)} \bigg(\frac{2}{3}\bigg) \to \bigg(\frac{1}{12}, \frac{1}{2}\bigg) \bigg(\frac{1}{8}, \frac{1}{8}, \frac{1}{2}\bigg) \bigg(\frac{1}{8}, \frac{3}{8}, \frac{1}{2}\bigg) \bigg(\frac{2}{3}\bigg)
\end{equation}
We underline the pile that the ``add and break'' operation is applied to at every step.

Theorem~\ref{thm:complete-is-optimal} implies that if one computes the sum of square roots of stick lengths in each pile, then all of these numbers cannot simultaneously exceed $1 + \sqrt{2}$, no matter how the player plays nor for how long.
Indeed, the lengths in each pile at the end are the squared magnitudes of the entries in each row of some normalization of $\bA^T$, where $T$ describes the branching structure of the player's choices; it is easily verified that any normalization can be obtained in this way by breaking sticks appropriately.

It is an intriguing question whether the construction of high-discrepancy matrices, no longer just unsatisfiable ones, may be encoded in some variant of such a game.

\subsection{Conjectural general dual certificate}
\label{sec:conj-general}

A final natural question, which we leave open, is whether our lower bound is optimal over $\whbA$ formed by normalizing \emph{any} $\bA^F$ for an unsatisfiable CNF formula $F$.
We conjecture that this is the case, and propose the following dual certificate construction to verify this claim, generalizing our use of read-once resolution proofs in the proof of Theorem~\ref{thm:complete-is-optimal}.

Let $T$ be a general resolution proof of unsatisfiability of $F$.
We use $T$ to define a dual certificate $\bm \mu$.
For a clause $C$ of $F$, let $\occ(C) \subseteq L(T)$ denote the set of leaves of $T$ that are labelled by $C$.
Define the \emph{effective depth} of $C$, denoted $d(C)$, by the relation
\begin{equation}
    2^{-d(C)} \colonequals \sum_{\ell \in \occ(C)} 2^{-\depth(\ell)}.
\end{equation}
Note that the effective depth obeys a coalescence property: the effective depth of a clause occurring twice at depth $k$ is $k - 1$, and similarly for further occurrences.
We then define, by analogy with $\bm\mu$ from the proof of Theorem~\ref{thm:complete-is-optimal},
\begin{align}
  \mu_{C}^{(0)} &\colonequals 2^{-d(C) / 2}, \\
  \mu_{C} &\colonequals \frac{\mu_{C}^{(0)}}{\sum_{C^{\prime} \in F} \mu_{C^{\prime}}^{(0)}}.
\end{align}
We then conjecture that the objective value of this choice of $\bm\mu$ is bounded by that achieved by the complete binary tree formula.
\begin{conjecture}
    \label{conj:dual-cert}
    For any unsatisfiable formula $F$ and resolution proof $T$, for $\bm\mu$ constructed as above,
    \begin{equation}
        \sqrt{\sum_{C \in F} \mu_C^2} + \sum_{z \in \var(F)} \sqrt{\sum_{C \in \partial z} \mu_C^2} = \frac{1 + \sum_{z \in \var(F)} \sqrt{\sum_{C \in \partial z} 2^{-d(C)}}}{\sum_{C \in F} 2^{-d(C) / 2}} < 1 + \sqrt{2}.
    \end{equation}
\end{conjecture}
\noindent
As mentioned earlier, this would also imply a positive solution of Problem~\ref{prob:unsat-opt} for unsatisfiable NAE-CNF formulas with constant $K_{\mathsf{unsat}} = 2 + \sqrt{2}$ (though the optimal constant could be smaller and indeed could be $1 + \sqrt{2}$).
We note that, by the coalescence property, the numerator above may be bounded by $\sum_{z \in \var(F)} 2^{-d(z)}$ for an effective depth of variables defined similarly over their occurrences when internal vertices are labelled with the variables resolved at those vertices.
That is, it would suffice to show the weaker bound
\begin{equation}
    \frac{1 + \sum_{z \in \var(F)} 2^{-d(z)/2}}{\sum_{C \in F} 2^{-d(C) / 2}} < 1 + \sqrt{2},
\end{equation}
which is directly analogous to the objective value considered in the proofs of Theorem~\ref{thm:complete-is-optimal}, with depth replaced by effective depth.
    
One potential approach to proving Conjecture~\ref{conj:dual-cert} is to define, following our second proof of Theorem~\ref{thm:complete-is-optimal}, a random map $f: F \into \var(F)$ from clauses to variables, such that we are able to control the difference between $d(C)$ and $d(f(C))$ on average.
However, it does not appear that a straightforward generalization of the partition into paths of $T$ from the proof of Theorem~\ref{thm:complete-is-optimal} furnishes a useful such map in this setting.

We note also that we do not need the claim to hold for \emph{any} resolution proof $T$, but rather just \emph{some} proof.
This may be helpful; for instance, \cite{Kullmann-2000-MatroidTheorySAT} describes how unsatisfiable $F$ having small deficiency admit resolution proofs consisting of several resolution proofs for $\MU(1)$ formulas ``glued'' together appropriately with a small number of further resolutions.
I have not yet been able to apply such results to generalize Theorem~\ref{thm:complete-is-optimal}, and numerical experiments and small examples appear to support the stronger conjecture stated above.

\section*{Acknowledgments}
\addcontentsline{toc}{section}{Acknowledgments}

I thank Afonso Bandeira, Pravesh Kothari, Cristopher Moore, and Daniel Spielman for helpful discussions during the preparation of this paper.

\appendix

\section{Second-order cone duality}
\label{app:socp}

In this section we justify the convex duality transformation of Lemma~\ref{lem:opt-norm}.
We take the approach of transforming the second-order cone program (SOCP) to a semidefinite program~(SDP) and computing the dual in that setting.
The key observation is that a constraint $\sqrt{\sum_{i = 1}^d y_i^2} \leq t$ on $\by \in \RR^d$ and with $t \geq 0$ is satisfied if and only if
\begin{equation}
    \bm M(\by, t) \colonequals \left[\begin{array}{cccc} t & y_1 & \cdots & y_d \\ y_1 & t & 0 & 0 \\ \vdots & 0 & \ddots & \vdots \\ y_d & 0 & \cdots & t\end{array}\right] \succeq \bm 0,
\end{equation}
as can be seen by taking the Schur complement with the bottom right $d \times d$ block.

In our case, we first rewrite the primal program in canonical form:
\begin{align}
  \delta(\whbA^F)
  &= \left\{ \begin{array}{ll} \text{maximize} & \min_{i = 1}^m \sum_{j \in C_i} a_{ij} \\ \text{subject to} & \sqrt{\sum_{i: j \in C_i} a_{ij}^2} \leq 1 \text{ for all } j \in [n] \end{array}\right\} \nonumber \\
  &= \left\{ \begin{array}{ll} \text{maximize} & c \\ \text{subject to} & c \leq \sum_{j \in C_i} a_{ij} \text{ for all } i \in [m], \nonumber \\
                                               & \sqrt{\sum_{i: j \in C_i} a_{ij}^2} \leq 1 \text{ for all } j \in [n] \end{array}\right\}
   \intertext{and then, letting $\ba_j \colonequals (a_{ij}: j \in C_i)$, we may rewrite as an SDP}
   &= \left\{ \begin{array}{ll} \text{maximize} & c \\ \text{subject to} & c \leq \sum_{j \in C_i} a_{ij} \text{ for all } i \in [m], \\
                                               & \bM(\ba_j, 1) \succeq \bm 0 \text{ for all } j \in [n] \end{array}\right\}.
\end{align}
We now introduce dual variables $\mu_i \geq 0$ for $i \in [m]$ for the linear constraints and $\bm\Lambda^{(j)} \succeq \bm 0$ with $\bm\Lambda^{(j)} \in \RR^{(\#\{i: j \in C_i\} + 1) \times (\#\{i: j \in C_i\} + 1)}$ for the semidefinite constraints.
The Lagrangian is
\begin{align}
  L(c, \ba; \bm\mu, \{\bm\Lambda^{(j)}\}_{j \in [n]})
  &\colonequals c - \sum_{i = 1}^m \mu_i\left(c - \sum_{j \in C_i} a_{ij}\right) + \sum_{j = 1}^n \left\langle \bM(\ba_j, 1), \bm \Lambda^{(j)} \right \rangle \nonumber \\
  &= \sum_{j = 1}^{n} \Tr(\bm \Lambda^{(j)}) + c\left(1 - \sum_{i = 1}^m \mu_i\right) + \sum_{i = 1}^m \sum_{j \in C_i} a_{ij}\left(\mu_i - 2\Lambda^{(j)}_{1,i + 1}\right).
\end{align}

Verifying the Slater conditions and writing the dual based on the rewritten Lagrangian above, we find the initial dual problem
\begin{equation}
    \label{eq:dual-1}
  \delta(\whbA^F)
  = \left\{\begin{array}{ll} \text{minimize} & \sum_{j = 1}^n \Tr(\bm\Lambda^{(j)}) \\ \text{subject to} & \bm \Lambda^{(j)} \in \RR^{(\#\{i: j \in C_i\} + 1) \times (\#\{i: j \in C_i\} + 1)}_{\succeq \bm 0}, \\ & \mu_{i} \geq 0, \\ & \Lambda^{(j)}_{1,i + 1} = \Lambda^{(j)}_{i + 1,1} = -\mu_{i} / 2 \text{ for all } i \in [m], j \in C_i, \\ & \sum_{i = 1}^m \mu_{i} = 1 \end{array} \right\}.
\end{equation}

We next observe that, for any fixed $\bm\mu$, we can optimize over all $\bm\Lambda^{(j)}$ explicitly.
Namely, this problem reduces to several independent problems of the form
\begin{equation}
    \begin{array}{ll} \text{minimize} & \Tr(\bm\Lambda) \\ \text{subject to} & \bm \Lambda \in \RR^{(d + 1) \times (d + 1)}_{\succeq \bm 0}, \\ & \Lambda_{1,i + 1} = \Lambda_{i + 1, 1} = x_i \text{ for all } i \in [d]. \end{array}
\end{equation}
Here, the positivity of $2 \times 2$ principal minors implies $\Lambda_{i + 1, i + 1} \geq x_{i}^2 / \Lambda_{1,1}$ for any feasible $\bm\Lambda$.
Moreover, for any $t > 0$ the choice of $\Lambda$ with $\Lambda_{1, 1} = t$, $\Lambda_{i + 1, i + 1} = x_i^2 / t$, $\Lambda_{1, i + 1} = \Lambda_{i + 1, 1} = x_i$, and all other entries equal to zero has $\bm\Lambda \succeq \bm 0$, since all principal minors are non-negative.
Therefore, the above is equivalent to minimizing over $t > 0$ the function $t + \sum_{i = 1}^d x_i^2 / t$, which is achieved at $t^{\star} = \sqrt{\sum_{i = 1}^d x_i^2}$ and has value $2t^{\star}$.

Applying this in \eqref{eq:dual-1}, we may eliminate the $\bm\Lambda^{(j)}$ variables, and we are left with the form of the dual given in Lemma~\ref{lem:opt-norm}:
\begin{equation}
    \delta(\whbA^F) = \left\{\begin{array}{ll} \text{minimize} & \sum_{j = 1}^n \sqrt{\sum_{i: j \in C_i} \mu_{i}^2} \\ \text{subject to} & \sum_{i = 1}^m \mu_{i} = 1 \end{array} \right\},
\end{equation}
where we note that the constraints $\mu_i \geq 0$ do not affect the value.

\addcontentsline{toc}{section}{References}
\bibliographystyle{alpha}
\bibliography{unsatisfiable-matrices}

\end{document}